\documentclass[11pt,a4paper]{article}
\usepackage[T1]{fontenc}
\usepackage{lmodern}
\usepackage[margin=1in]{geometry}
\usepackage{amsmath,amssymb,amsthm}
\usepackage{microtype}
\usepackage{enumitem}
\usepackage{cite}
\usepackage[hidelinks]{hyperref}
\hypersetup{pdfauthor={Kamran Fareed and Alina Lodhi},pdftitle={Understanding Loopless SVRG: Note V on Variance Reduction}}
\newtheorem{assumption}{Assumption}

\newtheorem{theorem}{Theorem}
\newcommand{\E}{\mathbb{E}}
\newcommand{\R}{\mathbb{R}}
\newcommand{\norm}[1]{\left\|#1\right\|}
\newcommand{\ip}[2]{\left\langle #1,#2\right\rangle}

\setlist[itemize]{leftmargin=1.4em,itemsep=2pt,topsep=4pt,parsep=0pt}
\title{Variance Reduction Note 5: The Convergence of Loopless SVRG}
\author{
{Kamran Fareed}$^{1}$ \quad
{Alina Lodhi}$^{1}$\\[0.5em]
\normalsize kamrafareed96@gmail.com, \quad xieyunqi27@gmail.com\\
\normalsize$^1$University of Agriculture Faisalabad
}
\date{}

\begin{document}
\maketitle
\vspace{-2em}
\begin{abstract}
A snapshot in loopless SVRG has no prescribed lifetime. At each step
the algorithm may keep using an old full gradient or compute a new one.
How does this choice enter the convergence proof, and how often is a
new full gradient worth its cost? For a strongly convex finite sum we
follow the age and objective gap of the snapshot. Its expected gap is
a geometrically weighted average of earlier iterate gaps. This identity
determines a weight for the snapshot in a one-step contraction argument
and gives an objective bound for the last iterate. A two-component
quadratic makes the estimator variance explicit, including on the step
immediately after a refresh. We also evaluate the work bound as a
function of the refresh probability and compare its minimizer with
the condition-number-independent choice $p=1/n$.
\end{abstract}

\section{The lifetime of a snapshot}
Consider $F(x)=n^{-1}\sum_{i=1}^n f_i(x)$. A full pass over the sum
at a point $w$ gives a vector $\nabla F(w)$ that can serve in several
successive updates. Classical SVRG \cite{svrg2013} organizes this
reuse in epochs. Note~4 of this series\footnote{Alina Lodhi,
\emph{Variance Reduction Note 4: The Convergence of SVRG},
\href{https://arxiv.org/abs/2611.05075v1}{arXiv:2611.05075v1}.}
analyzes an epoch of prescribed length and chooses the next snapshot
uniformly among its inner iterates. Its contraction compares the
objective values of two successive snapshots.

Kovalev, Horv\'ath, and Richt\'arik \cite{loopless2019} give the
snapshot a fresh chance of replacement at every iteration. There is
then no deterministic upper bound on the time between full passes when
$0<p<1$. A proof that required every reference point to be recent would
have to deal with these long intervals. We can instead ask what
information the snapshot carries after averaging over a single coin
flip. We will use the component-gradient estimate from Note~4,
Lemma~2, and develop the recursion for the random reference point here.

\paragraph{Which iterate is stored?}
Let $C_t$ be independent Bernoulli variables with success probability
$p\in(0,1]$. At initialization the snapshot is $w_0=x_0$.
At iteration $t$ a success stores $x_t$ as $w_{t+1}$, for use in
the next iteration. A failure leaves the stored point in place. To
distinguish the time of storage from the value of the vector, write
$w_t=x_{J_t}$ with $J_0=0$. Define the age $A_t=t-J_t$. It obeys
\[
 A_{t+1}=1+(1-C_t)A_t,\qquad A_0=0.
\]
In particular, $A_1=1$ under either outcome of the first toss. For
integers $1\leq k\leq t$, the event $A_t\geq k$ requires failures
in the last $k-1$ tosses and puts no conditions on earlier ones. Thus
\begin{equation}\label{eq:age}
 \Pr(A_t\geq k)=(1-p)^{k-1},\qquad
 \E A_t=\sum_{k=1}^{t}(1-p)^{k-1}
        =\frac{1-(1-p)^t}{p}.
\end{equation}
For $k=1$ the probability is one, including at $p=1$.
The mean age approaches $1/p$. Two snapshots of the same age can
nevertheless have quite different effects on the direction of a step,
depending on how far the intervening iterates have moved.

\section{The loopless update and its conditional variance}
\begin{assumption}\label{ass:convex}
Each $f_i:\R^d\to\R$ is differentiable and convex with an
$L$-Lipschitz gradient. The average $F$ is $\mu$-strongly convex,
where $0<\mu\leq L$. Let $x_\star$ be its minimizer,
$h(x)=F(x)-F(x_\star)$, and $\kappa=L/\mu$. Norms are Euclidean.
\end{assumption}

Fix a deterministic $x_0$ and compute $g_0=\nabla F(w_0)$.
At the start of step $t$ the algorithm has $x_t$, $w_t$, and the
cached vector $g_t=\nabla F(w_t)$. Draw a uniform index $i_t$
independently of the history and of the coin $C_t$. For a fixed
step size $\eta>0$, form
\begin{align}
 v_t&=g_t+\nabla f_{i_t}(x_t)-\nabla f_{i_t}(w_t),&
 x_{t+1}&=x_t-\eta v_t, \label{eq:update}\\
 (w_{t+1},g_{t+1})&=
 \begin{cases}
 (x_t,\nabla F(x_t)),& C_t=1,\\
 (w_t,g_t),& C_t=0 .
 \end{cases} \label{eq:refresh}
\end{align}
The output after $T$ steps is $x_T$. The reference in
\eqref{eq:refresh} is the point before the step, so a success at $t$
affects the direction starting at $t+1$. Keeping $w_t$ and $g_t$
requires $O(d)$ auxiliary storage. Write $\E_t$ for conditional
expectation given the history just before the two random draws.

\paragraph{Where the randomness enters.}
For fixed $x_t,w_t$, put $u_i=\nabla f_i(x_t)-\nabla f_i(w_t)$ and
$\bar u=n^{-1}\sum_i u_i$. Then $v_t-\nabla F(x_t)=u_{i_t}-\bar u$.
The estimator is conditionally unbiased for every pair of points, and
\begin{equation}\label{eq:localvariance}
 \E_t\norm{v_t-\nabla F(x_t)}^2
 =\frac1n\sum_{i=1}^n\norm{u_i}^2-\norm{\bar u}^2
 \leq L^2\norm{x_t-w_t}^2 .
\end{equation}
This identity separates the common change in the component gradients
from the part that varies with $i$. A large gradient at the reference
point can be harmless if the component gradients change by the same
amount on the way to $x_t$.

\subsection*{A quadratic on which the variance can be read off}
Take two equally likely components on the line. For $0\leq s<a$ and
$b\ne0$, define
\[
 f_+(x)=\tfrac12(a+s)x^2+bx,\qquad
 f_-(x)=\tfrac12(a-s)x^2-bx .
\]
Their average is $F(x)=ax^2/2$ and $x_\star=0$. We may take
$\mu=a$ and $L=a+s$. The linear terms give opposite nonzero
gradients at the solution. If $\xi_t$ is the sampled sign and
$d_t=x_t-w_t$, substitution into \eqref{eq:update} gives
\begin{equation}\label{eq:quad}
 v_t=ax_t+\xi_t s d_t,\qquad
 \E_t(v_t-ax_t)^2=s^2d_t^2 .
\end{equation}
The number $b$ cancels for every choice of $w_t$. When the
curvatures coincide ($s=0$), the method follows gradient descent
for any refresh schedule. When $s>0$, a reference left away from
the solution can produce a noisy step even at $x_t=0$.

The indexing of a refresh matters in this example. Given $C_t=1$,
we have $d_{t+1}=-\eta v_t$, so the variance in the next direction,
conditional on the state at $t+1$, is $s^2\eta^2v_t^2$. To see the
effect of the coin itself, start at $x_0=w_0=1$ with $a=2$, $s=1$
and $\eta=1/8$. The first step is exact, so $x_1=3/4$ and $w_1=1$.
If the next index selects the plus component, then $v_1=5/4$ and
$x_2=19/32$. A refresh now sets $w_2=3/4$, giving conditional
variance $(19/32-3/4)^2=25/1024$ at the next step. A failure keeps
$w_2=1$, giving $(19/32-1)^2=169/1024$. The two paths agree on
$x_2$. Their different reference points account for the entire
difference in variance.

\section{From snapshot survival to a last-iterate bound}
Let $S_t=h(w_t)$ and $D_t=\norm{x_t-x_\star}^2$. The coin acts
directly on $S_t$. Its two outcomes give
\begin{equation}\label{eq:memory}
 \E_t S_{t+1}=(1-p)S_t+p\,h(x_t).
\end{equation}
Taking expectations and expanding this recursion through $t$ steps
yields
\begin{equation}\label{eq:history}
 \E S_t=(1-p)^t h(x_0)
       +p\sum_{j=0}^{t-1}(1-p)^{t-1-j}\E h(x_j).
\end{equation}
The weights sum to one. For example, an objective gap incurred $k$
steps earlier enters the sum with weight $p(1-p)^{k-1}$.
Equation~\eqref{eq:history} follows by conditioning one step at a
time. The future trajectory depends on which refreshes occurred,
so conditioning on the entire trajectory would change the problem.

\paragraph{The estimate inherited from classical SVRG.}
Conditioning fixes the two arguments of the estimator in
\eqref{eq:update}. Under Assumption~\ref{ass:convex}, Lemma~2 of
Note~4 therefore applies at this step and gives
\begin{equation}\label{eq:import}
 \E_t\norm{v_t}^2\leq4L\bigl[h(x_t)+S_t\bigr].
\end{equation}
This second-moment bound can be used directly in a squared-distance
update. In contrast, the exact variance identity
\eqref{eq:localvariance} describes what happens at a particular
pair of iterates, as in \eqref{eq:quad}. Strong convexity gives
$\ip{\nabla F(x_t)}{x_t-x_\star}\geq h(x_t)+\mu D_t/2$.
Expanding $D_{t+1}$, using $\E_t v_t=\nabla F(x_t)$, and applying
\eqref{eq:import} now yields
\begin{equation}\label{eq:dist}
 \E_t D_{t+1}\leq (1-\mu\eta)D_t
 -(2\eta-4L\eta^2)h(x_t)+4L\eta^2 S_t .
\end{equation}
The last term allows a poor reference to offset the progress made
by the current direction.

\paragraph{Choosing the weight of a stale snapshot.}
Consider $\Phi_t=D_t+cS_t$ with $c>0$. The snapshot term in
\eqref{eq:dist} costs $4L\eta^2 S_t$ at every step. From
\eqref{eq:memory} we recover $cpS_t$ in expected decrease, at the
price of $cp\,h(x_t)$. To retain a contraction margin of $p/2$ on
$cS_t$, we can set $cp/2=4L\eta^2$. The coefficient of $h(x_t)$
in the combined bound then becomes $-2\eta+12L\eta^2$.
In particular $\eta=1/(8L)$ leaves a decrease of $\eta h(x_t)/2$.
This calculation gives a single weight that works for every $p$ once
it is scaled by $1/p$.

\begin{theorem}[Contraction of the iterate and reference]\label{thm:linear}
Suppose Assumption~\ref{ass:convex} holds. Set
\[
 \eta=\frac1{8L},\qquad c=\frac{8L\eta^2}{p},\qquad
 \delta=\min\{\mu\eta,p/2\}.
\]
Then $\E_t\Phi_{t+1}\leq(1-\delta)\Phi_t$. In particular, for any
integer $T\geq0$ the returned iterate satisfies
\begin{equation}\label{eq:last}
 \E h(x_T)\leq \tfrac L2 (1-\delta)^T\Phi_0 .
\end{equation}
\end{theorem}
\begin{proof}
Adding \eqref{eq:dist} to $c$ times \eqref{eq:memory} with the stated
parameters gives
\[
 \E_t\Phi_{t+1}\leq(1-\mu\eta)D_t+(1-p/2)cS_t
                         -\tfrac{\eta}{2}h(x_t).
\]
Both $D_t$ and $S_t$ are nonnegative, and $h(x_t)\geq0$, which
proves the conditional bound. Iterated expectation gives
$\E\Phi_T\leq(1-\delta)^T\Phi_0$. Since $F$ is $L$-smooth and
$\nabla F(x_\star)=0$, $h(x_T)\leq LD_T/2\leq L\Phi_T/2$.
\end{proof}

\subsection*{How much work should be spent on refreshes?}
With two component evaluations per direction and a fresh full pass
at every success, the number of evaluations at a fixed horizon is
\begin{equation}\label{eq:cost}
 N_T=n+2T+nK_T,\qquad K_T=\sum_{t=0}^{T-1}C_t,\qquad
 \E N_T=n+(2+np)T .
\end{equation}
Here $K_T$ has the binomial distribution with parameters $T$ and $p$.
The count includes the initial full pass and allows a refresh on the
final step. Reusing coincident evaluations or omitting that final
refresh can only reduce the count.

For a relative potential target $\E\Phi_T\leq r\Phi_0$, $0<r<1$,
the theorem permits $T=\lceil\max\{8\kappa,2/p\}\log(1/r)\rceil$.
Before rounding, the coefficient of $\log(1/r)$ in the expected
work is
\begin{equation}\label{eq:work}
 W(p)=(2+np)\max\{8\kappa,2/p\}
 =\begin{cases}
       4/p+2n, & 0<p\leq 1/(4\kappa),\\
       16\kappa+8\kappa np, & 1/(4\kappa)\leq p\leq1 .
   \end{cases}
\end{equation}
On the first interval, less frequent full passes increase the number
of iterations needed to replace an old reference. On the second,
the iteration bound already depends only on $\kappa$, so a larger
$p$ raises the work at each iteration. The minimum of this particular
upper-bound coefficient occurs at $p_\star=1/(4\kappa)$ and equals
$16\kappa+2n$. The minimizing value refers to the work model
\eqref{eq:cost} and the constants in Theorem~\ref{thm:linear}.

When $\mu$ is unknown, $p=1/n$ gives $W(1/n)=3\max\{8\kappa,2n\}$.
For $n\geq4\kappa$ this is $6n$, at most $3W(p_\star)$.
For $n\leq4\kappa$ it is $24\kappa$, at most
$\tfrac32 W(p_\star)$. Thus the simple rule stays within a constant
factor on both sides of the transition. These comparisons concern
relative reduction of the potential. An absolute objective target
also depends on its starting value, which changes with $p$.

\paragraph{Conversion to objective accuracy.}
At $p=1/n$ we have $c=n/(8L)$. Write $\Delta_0=h(x_0)$.
Since $w_0=x_0$ and $D_0\leq2\Delta_0/\mu$,
\[
 \tfrac L2\Phi_0
   =\tfrac L2D_0+\tfrac n{16}\Delta_0
   \leq(\kappa+n/16)\Delta_0 .
\]
For $0<\tau<\Delta_0$, \eqref{eq:last} and
$(1-\delta)^T\leq e^{-\delta T}$ give
\begin{equation}\label{eq:target}
 T=\left\lceil\max\{8\kappa,2n\}
      \log\frac{(\kappa+n/16)\Delta_0}{\tau}\right\rceil
 \quad\Longrightarrow\quad \E h(x_T)\leq\tau,\qquad
 \E N_T=n+3T
\end{equation}
for the implementation counted in \eqref{eq:cost}. If
$\Delta_0\leq\tau$, the initial point suffices. Running the updates
at $p=1/n$ requires $n$ and a smoothness bound for setting $\eta$.
Choosing a horizon from \eqref{eq:target} additionally uses a lower
bound on the strong convexity constant and an upper bound on the
initial objective gap.

\subsection*{Choosing a horizon under a fixed gradient budget}
Suppose that only $B$ component-gradient evaluations are available.
At $p=1/n$, setting $n+3T=B$ spends $B$ evaluations in expectation.
The probability of finishing those $T$ steps depends on how many full
passes can actually fit. In what follows we count evaluations as in
\eqref{eq:cost}, including the last refresh, and take $T$ to be fixed
before drawing any indices or coins.

After reserving $n$ evaluations for initialization and $2T$ for the
directions, the number of full passes still affordable is
\[
 m_B(T)=\left\lfloor\frac{B-n-2T}{n}\right\rfloor .
\]
For $0<p<1$ the probability of completing the run within budget is
\begin{equation}\label{eq:budgetprob}
 q_B(T):=\Pr(N_T\leq B)
       =\sum_{j=0}^{\min\{T,m_B(T)\}}
          {T\choose j}p^j(1-p)^{T-j}.
\end{equation}
Indeed $N_T\leq B$ precisely when $K_T\leq m_B(T)$. An empty sum is
zero, and $m_B(T)\geq T$ makes the probability one. At $p=1$ the
cost is deterministic and the completion condition is simply
$n+(n+2)T\leq B$.

Given an acceptable overrun probability $\alpha\in(0,1)$ and $B\geq n$,
choose the largest integer $T$ with $q_B(T)\geq1-\alpha$. This takes
only a one-dimensional search up to $\lfloor(B-n)/2\rfloor$.
To see why the test is monotone, generate one sequence of refresh
coins and use its first $T$ tosses for each candidate horizon. On
every such sequence $N_{T+1}-N_T=2+nC_T>0$, so a failure to fit at
$T$ is also a failure to fit at $T+1$.

\paragraph{One more step can cost an entire refresh slot.}
Take $n=100$, $p=1/100$ and $B=600$. Equation~\eqref{eq:budgetprob}
gives the following values, rounded to four decimal places:
\begin{center}
\begin{tabular}{rccc}
\hline
Steps $T$ & Mean cost $\E N_T$ & Affordable refreshes $m_B(T)$
          & Completion probability $q_B(T)$\\
\hline
99  & 397 & 3 & 0.9822\\
100 & 400 & 3 & 0.9816\\
101 & 403 & 2 & 0.9188\\
\hline
\end{tabular}
\end{center}
At $T=100$, initialization and the directions leave exactly $300$
evaluations for refreshes. At $T=101$ they leave $298$, enough for
only two passes. Thus $100$ is the largest horizon that fits with
probability at least $0.95$, even though the mean cost at $101$ is
well below $600$.

\paragraph{Relating the budget to objective accuracy.}
A large completion probability alone says little about the typical
objective gap on completed runs. A refresh changes both the cost and
the future iterates. One way to state both guarantees without assuming
independence is to bound the probability of either failure. Let
$M_T=(L/2)(1-\delta)^T\Phi_0$, the expectation bound in
\eqref{eq:last}. For any target $\tau>0$,
\begin{equation}\label{eq:joint}
 \Pr\{N_T\leq B,\ h(x_T)\leq\tau\}
 \geq 1-\bigl[1-q_B(T)\bigr]-M_T/\tau .
\end{equation}
The first failure event is exceeding the budget. For the second,
Markov's inequality gives $\Pr\{h(x_T)>\tau\}\leq M_T/\tau$ because
$h$ is nonnegative. Adding the two failure probabilities proves
\eqref{eq:joint}. For instance, if $q_B(T)\geq1-\alpha$ and
$M_T\leq\beta\tau$, a run completes within budget and meets the
objective target with probability at least $1-\alpha-\beta$.
This separates two practical choices: $n$, $p$, and $B$ determine
which horizons fit, while $\mu$, $L$, and the starting potential
determine which of those horizons certify the desired accuracy.

\section{Recommended reading}

\paragraph{Variance-reduction methods.}
Basic guarantees \cite{ghadimi2013,gower2020,unifiedsgd2020}; gradient memories and control
variates \cite{sag2012,sdca2013,svrg2013,s2gd2013,saga2014,finito2014,proxsvrg2014}; sampling and acceleration \cite{catalyst2015,katyusha2016,svrgpp2016,katyushax2018,sagasampling2019,jacsketch2018,svrgsd2018,loopless2019,looplesssampling2021}; nonconvex and recursive
estimators \cite{reddi2016,allenzhu2016,natasha2017,sarah2017,isarah2018,scsg2017,spider2018,snvrg2018,spiderboost2019,proxsarah2020,zerosarah2021,page2021,nonconvexsampling2019}; momentum and adaptive schedules \cite{storm2019,hybrid2019,stormplus2021,superadam2021,metastorm2022,adaspider2022,jiang_adaptive2024,adasvrg2021,almosttune2020,adaptiveextra2022,alfpg2026};
second-order extensions \cite{ssrgd2019,natasha22018,svrc2019,svrcnewton2018}; lower bounds \cite{lower2023}.

\paragraph{Compositional optimization.}
Two-level  approximation \cite{scgd2014,asc2017,nasa2020,nearlyeasy2020}; finite-sum and recursive
estimators \cite{lian2017,compavg2018,civr2019,sarahcomp2019}; multilevel tracking \cite{nested2019,jiang_multilevel2022,levelindependent2020,jiang_revisiting2025};
coupled and partially observed blocks \cite{fccotheory2022,jiang_msvr2022,weakcomp2023,jiang_fcco2026,momentumcomp2025}; projection-free
methods \cite{vrfw2016,projectioncomp2022,jiang_projection2024,scgpp2019,jiang_projection2026}; statistical-model applications and federated
extensions \cite{jiang_models2023,fednest2022,jiang_models2026,localcomp2022}.

\paragraph{Sign-based optimization.}
Quantization, signSGD, majority voting, and error feedback \cite{qsgd2017,terngrad2017,signsgd2018,majority2019,errorfeedback2019};
variance-reduced sign updates and Lion \cite{jiang_sign2024,lion2023,jiang_lion2026}; compressed gradient
tracking \cite{diana2019,marina2021,jiang_distributed2026,ef212021,dasha2022}. Related distributed control-variate methods include
SCAFFOLD, FedPAGE, decentralized hybrid estimators, and ProxSkip \cite{scaffold2020,fedpage2021,hybriddecentral2021,proxskip2022}.

\paragraph{Minimax optimization.}
Finite-sum saddle-point methods, matrix games, and variance-reduced
extragradient \cite{saddle2016,matrixgames2019,svre2019,catalystminimax2020}; recursive, hybrid, adaptive, and decentralized
nonconvex minimax methods \cite{sreda2020,hybridminimax2020,adgda2023,dream2024}; a variational-inequality treatment \cite{vi2022}.


\end{document}